\documentclass[11pt]{article}
\usepackage{amssymb}
\usepackage{amsmath}
\usepackage{amsthm}
\usepackage{enumerate}
\usepackage{color}
\newtheorem{def1}{Definition}[section]
\newtheorem{th1}{Theorem}[section]

\newtheorem{conc1}{Conclusion}[section]

\newtheorem{ex1}{Example}[section]

\numberwithin{equation}{section}

\newcommand{\Zpcz}[2]{{#1}\! \! \! \! \!  _{_{_{_{#2}}}} }   
\newcommand{\dps}{\displaystyle}
\newcommand{\pcz}[1]{\frac{\partial}{\partial{#1}}}
\newcommand{\Pcz}[2]{\frac{\partial{#1}}{\partial{#2}}}
\newcommand{\Ppcz}[2]{\frac{\partial^2{#1}}{\partial{#2}^2}}
\newcommand{\PPcz}[3]{\frac{\partial^2{#1}}{\partial{#2}\partial{#3}}}
\newcommand{\p}[1]{\partial _{#1}}
\newcommand{\klam}[1]{\left\{ \begin{array}{ll}  #1  \end{array} \right.}
\newcommand{\mb}[1]{\mathbb{#1}}
\newcommand{\wt}{\widetilde}

\newcommand{\ve}{\varepsilon}
\newcommand{\vp}{\varphi}

\begin{document}

\title{ Equivalence and integrability of second order linear ODEs}
\author{Ivan Tsyfra$^{(1)}$, \ Tomasz Czy\.{z}ycki$^{(2)}$}

\maketitle

\noindent $^{(1),(2)}$ { Institute of Mathematics, University of Bia\l ystok, Poland, \\
 e-mail: $^{(1)}$tsyfra@math.uwb.edu.pl, \ $^{(2)}$tomczyz@math.uwb.edu.pl } \\[4pt]


\vspace{1cm}

\begin{abstract}
We consider a class of linear ODEs of  second order with variable coefficients and construct its Lie algebra of  Lie group of equivalence transformations. Further we find invariants and differential invariants of this Lie algebra and by using them we formulate criteria of equivalence of the  equations under consideration.
These criteria enable us to characterize some classes integrable in quadratures.
\end{abstract}

\section{Introduction}

The notion of equivalence transformations is defined for a class of differential equations. It enable one to investigate properties of this class and find the change of variables that transforms one equation to another.
The equivalence transformations have already been originally used by S. Lie for group classification of ordinary second order  differential equations in \cite{L1}.
The invariants of equivalence transformations and infinitesimal technique of generating these invariants were introduced by L. Ovsiannikov in \cite{O1, O2}.

In this paper we find the Lie algebra of the Lie group of equivalence transformations of the second order ODEs.
Using this algebra and its invariants we formulate criteria of equivalence and integrability of these equations.
This new approach is based on idea of invariant combinations of semi-invariants and is discovered by  I. Tsyfra.(see Theorem 4.1).
It allows us to avoid difficult calculations and to compare only few simple expressions to establish whether two equations are equivalent or not. If the equivalence conditions hold, it is possible to construct exact transformations which reduce one equation to another.

In Section 2, using the infinitesimal technique of constructing the Lie algebra of a given Lie group, we find the Lie group of equivalence transformations of the class of linear second order ODEs.

In Section 3, following the procedure based on  Tresse theorem, we construct invariants, differential invariants and invariant differentiation operators for each subgroup of the  Lie group of equivalence transformations constructed in Section 2.

In Section 4,  by the usage of invariants obtained the criteria of equivalence of the studied equations is given.

In Section 5,  criteria of integrability is applied to a certain class of equations.

\section{Equivalence transformations of the second order linear ODEs }

Consider the family of linear ODEs of the second order, parametrized by two functions $a(x), b(x)$ in the form:
\begin{equation}\label{Eq}
y''+a(x)y'+b(x)y=0.
\end{equation}
We have following
\begin{def1}
Equivalence transformation of the family of second order linear ODEs  (\ref{Eq}) is a regular change of variables
\begin{equation}
\wt{x}=\alpha(x,y), \ \wt{y}=\beta(x,y),
\end{equation}
preserving the set
\begin{equation} \label{SetEq}
\Omega=\{y''+a(x)y'+b(x)y=0: \ y,a,b:\mb{R}\to\mb{R} \}.
\end{equation}
\end{def1}

Denote by $G$ the Lie group of equivalence transformations of $\Omega$.
We are searching for a Lie algebra ${\cal A}$ of $G$.
Treating $a,b$ as new dependent variables we are searching for the symmetry operators
\begin{equation}
X=\xi(x,y)\p{x}+\eta(x,y)\p{y}+\mu^1(x,y,a,b)\p{a}+\mu^2(x,y,a,b)\p{b}
\end{equation}
 of the system of equations
\begin{equation} \label{ukl1}
\begin{cases} y''+ay'+by=0. \\[-9pt]
\cr \dps{\Pcz{a}{y}=0, \ \Pcz{b}{y}=0}. \\[-7pt]
\cr \dps{\PPcz{a}{x}{y}=0, \ \Ppcz{a}{y}=0, \ \PPcz{b}{x}{y}=0, \ \Ppcz{b}{y}=0,}
\end{cases}
\end{equation}
We calculate the first and second order extension of operator $X$ to the space with coordinates\\
$(x,y,y',y'',a,a_x,a_y,b,b_x,b_y,a_{xx},a_{xy},a_{yy},b_{xx},b_{xy},b_{yy} )$:
\begin{equation}
\Zpcz{X}{1}=X+\zeta^1\p{y'}+\nu^1_1\p{a_x}+\nu^1_2\p{a_y}+\nu^2_1\p{b_x}+\nu^2_2\p{b_y}+\nu^3_1\p{c_x}+\nu^3_2\p{c_y},
\end{equation}
\begin{equation}
\Zpcz{X}{2}=\Zpcz{X}{1}+\zeta^{11}\p{y''}+\nu^1_{11}\p{a_{xx}}+\nu^1_{12}\p{a_{xy}}+\nu^1_{22}\p{a_{yy}}+ \nu^2_{11}\p{b_{xx}}+\nu^2_{12}\p{b_{xy}}+\nu^2_{22}\p{b_{yy}},
\end{equation}
where
$$\zeta^1=D_x(\eta)-y'D_x(\xi), $$
$$\zeta^{11}=D_x(\zeta^1)-y''D_x(\xi), $$

$$\nu^1_1=\wt{D}_x(\mu^1)-a_x\wt{D}_x(\xi)-a_y\wt{D}_x(\eta), $$
$$\nu^2_1=\wt{D}_x(\mu^2)-b_x\wt{D}_x(\xi)-b_y\wt{D}_x(\eta), $$
$$\nu^1_2=\wt{D}_y(\mu^1)-a_x\wt{D}_y(\xi)-a_y\wt{D}_y(\eta), $$
$$\nu^2_2=\wt{D}_y(\mu^2)-b_x\wt{D}_y(\xi)-b_y\wt{D}_y(\eta), $$
$$\nu^1_{11}=\wt{D}_x(\nu^1)-a_{xx}\wt{D}_x(\xi)-a_{xy}\wt{D}_x(\eta),\ etc. $$

and total derivatives are in the form
$$D_x=\pcz{x}+y'\pcz{y} $$
$$\wt{D}_x=\pcz{x}+a_x\pcz{a}+b_x\pcz{b}+c_x\pcz{c} $$
$$\wt{D}_y=\pcz{y}+a_y\pcz{a}+b_y\pcz{b}+c_y\pcz{c}. $$
Taking into account equations (\ref{ukl1}), we have $\dps{\wt{D}_y=\pcz{y}} $.

 Acting by $\Zpcz{X}{2}$ \ on system (\ref{ukl1}) we obtain the determining equations \begin{equation}
\klam{ \zeta^{11}+\mu^1 y'+a\zeta^1+\mu^2y+b\eta \ \big|_{(\ref{ukl1})}=0 \\[6pt]
 \nu^1_2=\nu^2_2=\nu^1_{12}=\nu^1_{22}=\nu^2_{12}=\nu^2_{22}\ \big|_{(\ref{ukl1})}=0 }
\end{equation}
and
$$\mu^1_y=\mu^2_y=0 \quad \wedge \quad \xi_y=0. $$

The first equation has the form:
\begin{equation} \label{determ}
\eta_{xx}+a\eta_x+b\eta-by\eta_y+2by\xi_x+\mu^2y+y'(2\eta_{xy}-\xi_{xx}+a\xi_x+\mu^1)+y'\,^2\eta_{yy}  \ \big|_{(\ref{ukl1})}=0.
\end{equation}

Now we separate variables with respect to $y'$ and after integration of separated system we obtain
\begin{th1}
The Lie algebra ${\cal A}$ of the Lie group $G$ is generated by the operators:
\begin{equation} \label{oper1}
X=\xi(x)\p{x}+A(x)y\p{y}+(\xi''-a\xi'-2A')\p{a}-(A''+aA'+2b\xi')\p{b}
\end{equation}
 where $\xi(x), A(x)$ are arbitrary smooth functions.
 \end{th1}

\noindent Now we describe two subgroups of $G$, generated by operators with one arbitrary function \\[4pt]
1) $\xi(x)$:
\begin{equation}\label{X1}
X_1=\xi(x)\p{x}+(\xi''-a\xi')\p{a}-2b\xi'\p{b}
\end{equation}
Equivalence group transformations is of the form:
 \begin{equation}\label{zam1}
 \begin{cases}
\wt{x}=\wt{\xi}^{-1}(\ve+\wt{\xi}(x))=\alpha(x) \cr \wt{y}=y,
\end{cases}
\end{equation}
where $\dps{\wt{\xi}(s)=\int\frac{ds}{\xi(s)}}$. \\[4pt]

\noindent 2) $A(x)$:
\begin{equation}\label{X2}
X_2=A(x)y\p{y}-2A'\p{a}-(A''+aA')\p{b}
\end{equation}
Equivalence group of transformations is of the form:
 \begin{equation}\label{zam2}
 \begin{cases}
\wt{x}=x \cr
\wt{y}=ye^{\ve A(x)}
\end{cases}
\end{equation}
Here $\ve$ is a group parameter.

\section{Invariants of equivalence transformations of the family of second order linear ODEs}

In this section we calculate invariants and differential invariants of the Lie group of equivalence transformations constructed in Section 2. It occures that operator (\ref{oper1}) has no invariants.
This means that invariants for the entire Lie group $G$ do not exist.
Therefore we find invariants of subgroups of $G$.
Further by using Tresse theorem \cite{Tr, Ib1} we construct basis of differential invariants and operators of invariant differentiation for each subgroup. \\
Let us consider subgroups generated by operators (\ref{X1}), (\ref{X2}) and look for expressions satisfying the condition:
$$\Zpcz{X}{m}\omega(x,y,y',a,b,a',b',...,a^{(m)},b^{(m)})=0. $$
The basic results are the following: \\[4pt]
for $X_1$:
$$ \omega_{11}=y, \quad \omega_{12}=b\xi^2, \quad \omega_{13}=\dps{\frac{a^2}{b}+\frac{ab\,'}{b^2}+\frac{b'\,^2}{4b^3} } $$
Operator of invariant differentiation has the form \ $\dps{Q_1=\frac{1}{\sqrt{|b|}}D_x}$ \\[4pt]
for $X_2$:
$$\omega_{21}=x, \quad \omega_{22}=\dps{a+2\frac{A'}{A}\ln y'}, \quad \omega_{23}=a^2-4b+2a' $$
The operator of invariant differentiation is \ $Q_2=D_x$ \\[4pt]
Hence by invariant differentiation one can obtain infinite sequence of differential ivariants for these subgroups.
E.g. for $X_1$ we have
$$\omega_{14}=Q_1(\omega_{13})=\dps{\frac{1}{\sqrt{|b|}}D_x(\omega_{13})}, \quad \omega_{15}=Q_1(\omega_{14})
\quad {\rm etc.} $$

\section{Criteria of equivalence of the second order linear ODEs }

In this section we formulate criteria of equivalence of two equations:
\begin{equation}\label{eq1}
y''(x)+a_1(x)y'(x)+b_1(x)y(x)=0
\end{equation}
 and
 \begin{equation}\label{eq2}
 z''(t)+a_2(t)z'(t)+b_2(t)z(t)=0.
 \end{equation}
 We say that equations (\ref{eq1}), (\ref{eq2}) are equivalent if there exists regular change of variables
 $$\klam{x=\vp(t,z) \\ y=\psi(t,z)} $$
 which transforms one equation to another. \\
These criteria are based on the following theorem which is due to I.Tsyfra and the proof for general case of partial differential equations will be presented in his separate paper under preparation.

\begin{th1} \label{thniezmkonstr}
Let \ $\omega_1, \omega_2,...,\omega_n$, where $\omega_i=\omega_i(x,a,b,a',b',...,a^{(k)},b^{(k)} )$ for a certain $k\in\mathbb{N}$, \ are the differential invariants of the Lie group with generator (\ref{oper1}).  \\
Assume that the equation (\ref{eq1}) \ with $a=a_1(x), b=b_1(x)$ can be transformed to the equation (\ref{eq2}) with $a=a_2(t), b=b_2(t)$ by the group of equivalence transformations.
Then for any smooth function  $H$ satisfying the condition
\begin{equation} \label{niezmkomb11}
 H(\omega_1, \omega_2,...,\omega_n)\big|_{ {a=a_1(x)}\atop{b=b_1(x)} } =\lambda=const.
 \end{equation}
it follows that
\begin{equation} \label{dod}
 H(\omega_1, \omega_2,...,\omega_n)\big|_{ {a=a_2(t)}\atop{b=b_2(t)} } =\lambda.
\end{equation}
\end{th1}

Sufficient conditions of equivalence for such two equations are also expressed in terms of (\ref{niezmkomb11}) and (\ref{dod}).

Generally speaking, if there are at least two functionally independent differential invariants $\omega_1, \omega_2$, then there exists a function $H$,
satisfying  condition (\ref{niezmkomb11}). \\

\begin{ex1} {\rm}
Let us  consider two equations \ $\dps{y''(x)+\frac{1}{x}y'(x)+4xy(x)=0}$ \ and \ $y''(t)+4e^{-3t}y(t)=0$. \\[3pt]
The question is, are these equations equivalent and, if so, what are the equivalence transformations beetwen them.
The direct way to verify this fact is complicated, but we can use Theorem \ref{thniezmkonstr} and compare invariants and invariants combinations of the Lie group of equivalence transformations for these equations. \\
The invariants of $X_1$ connected with the first equation are the following
$$\omega_{13}=\frac{9}{16x^3}, \quad \omega_{14}=\frac{27}{32x^4\sqrt{x}} $$
Invariants connected with the second one are
$$\omega_{23}=\frac{9}{16}e^{3t}, \quad \omega_{14}=\frac{27}{32}e^{\frac{9}{2}t}  $$
Hence we obtain invariant combinations
$$H_1(\omega_{13},\omega_{14})=\frac{4\omega_{13}^3}{\omega_{14}^2}=1 $$
$$H_2(\omega_{23},\omega_{24})=\frac{4\omega_{23}^3}{\omega_{24}^2}=1 $$
and conclude that given equations are equivalent with respect to transformations generated by $X_1$.
We can also calculate these transformations by comparing $\omega_{13}$ and $\omega_{23}$ which give
$x=e^{-t}$. It is easy to verify that this transformation with invariant $y$ is the  equivalence transformation searched for.
\end{ex1}

\section{Application to integrability of the linear second order ODEs }

Investigations of the problem of integrability of the linear second order ODEs are based on the following observation:
Assume that equation (\ref{eq1}) is integrable and equation (\ref{eq2})
is equivalent to (\ref{eq1}) by means of explicit equivalence transformations.
Then equation (\ref{eq2}) is integrable as well.

\noindent In this context a question arises: When the general second order linear ODE (\ref{Eq}) can be transformed into the form
\begin{equation} \label{eq3}
y''+V(x)y=0.
\end{equation}

Theorem \ref{thniezmkonstr} implies

\begin{conc1} \label{thniezmkonstr2}
Let \ $\omega_1, \omega_2,...,\omega_n$ be the differential invariants of the Lie group with generator (\ref{oper1}).  \\
Assume that  equation (\ref{eq1}) \ with $a=a_1(x), b=b_1(x)$ can be transformed to  equation (\ref{eq3}) with $V=V_1(x)$ by the group of equivalence transformations. Then for any smooth function $H$ satisfying the condition
\begin{equation} \label{niezmkomb111}
 H(\omega_1, \omega_2,...,\omega_n)\big|_{ {a=a_1(x)}\atop{b=b_1(x)} } =\lambda=const
 \end{equation}
 one has
$$ H(\omega_1, \omega_2,...,\omega_n)\big|_{ {a=0}\atop{b=V(x)} } =\lambda.$$
\end{conc1}

To answer the additional question whether the equation $y''+V_1(x)y=0 $ can be transformed into equation $y''+V_2(x)y=0 $ one has to use the group of equivalence transformations with the infinitesimal generator
\begin{equation} \label{oper3}
X=2\xi(x)\p{x}+\xi'y\p{y}-(\xi'''+4\xi'V)\p{V}.
\end{equation}
For the group of transformations generated by (\ref{oper3}) the invariant differentiation operator has the form $Q=\xi(x)D_x$ and differential invariants are:
$$\omega_1=V\xi^2+\frac{1}{2}\xi\xi''-\frac{1}{4}\xi'^2, \quad \omega_2=Q\omega_1. $$
In this case differential invariants depend on $x$ and therefore Theorem \ref{thniezmkonstr} holds for the transformations generated by (\ref{oper3}) for a given $\xi(x)$ only, but not for arbitrary $\xi(x)$. \\[6pt]
As far as $\xi(x), A(x)$ are arbitrary functions, one can use the group transformations to produce the class of potentials $V_{\xi}(x), V_{A}(x)$ from a given potential $V$. If equation (\ref{eq3}) is integrable (nonintegrable) for this potential $V(x)$, then it is also integrable (nonintegrable) for any potential belonging to this class. If there exists a function $H$ such that:
\begin{equation} \label{niezmkomb1}
 H(\omega_1, \omega_2,...,\omega_n)\big|_{ V=V_1(x) } =\lambda_1=const.,
 \end{equation}
 for given $\xi(x)$ and $H(\omega_1, \omega_2,...,\omega_n)\big|_{ V=V_2(x) } \ne\lambda_1,$ then $V_1,V_2$ do not belong to the same class as follows from theorem \ref{thniezmkonstr}. \\[5pt]

\end{document}